\documentclass[a4paper,10pt]{article}
\usepackage{geometry}
\usepackage{cite}
\usepackage{breqn}
\usepackage[linesnumbered,ruled,vlined]{algorithm2e}
\usepackage{todonotes}
\usepackage{authblk}

\SetCommentSty{mycommfont}

\SetKwInput{KwInput}{Input}                
\SetKwInput{KwOutput}{Output}                

\newcommand{\figref}[1]{Fig.~\ref{fig:#1}}

\newcommand{\eqnref}[1]{Eq.~\eqref{eq:#1}}
\newcommand{\secref}[1]{Sec.~\ref{sec:#1}}

\newcommand{\algref}[1]{Algorithm~\ref{alg:#1}}
\newcommand{\eq}[1]{\eqref{eq:#1}}

\newcommand{\bvec}[1]{\boldsymbol{#1}}
\newcommand{\mat}[1]{\boldsymbol{#1}}

\newcommand{\diff}{{\rm d}}

\newcommand{\zL}{z_{L}}
\newcommand{\zR}{z_{R}}

\newcommand{\um}{\mu{\rm m}}

\title{VarRCWA: An Adaptive High-Order Rigorous Coupled Wave Analysis Method}
\author[1]{Ziwei Zhu}
\author[1]{Changxi Zheng\thanks{cxz@cs.columbia.edu}}
\affil[1]{Department of Computer Science, Columbia University, New York, New York 10027, USA}
\begin{document}
\maketitle
\begin{abstract}
    Semi-analytical methods, such as rigorous coupled wave analysis, have been
    pivotal for numerical analysis of photonic structures. 
    In comparison to other methods, they offer much faster computation, 
    especially for structures with constant cross-sectional shapes (such as metasurface units).
    However, when the cross-sectional shape varies even mildly (such as a taper), 
    existing semi-analytical methods 
    suffer from high computational cost.
    We show that the existing methods can be viewed as a zeroth-order approximation with respect to the structure's 
    cross-sectional variation. We instead derive a high-order perturbative expansion with respect to the cross-sectional
    variation. Based on this expansion, we propose a new semi-analytical method
    that is fast to compute even in presence of large cross-sectional shape variation.
    Furthermore, we design an algorithm that automatically discretizes the
    structure in a way that achieves a user specified accuracy level while at
    the same time reducing the computational cost.
\end{abstract}
\section{Introduction}





Numerical simulation is a fundamental tool for understanding photonic
structures. Among many popular methods, semi-analytical methods, such as
rigorous coupled wave analysis (RCWA)~\cite{moharam1981rigorous}, have been
widely used for analyzing such devices as
metasurfaces~\cite{divitt2019ultrafast}, gratings~\cite{mohamad2020fast} and
waveguides~\cite{zhu2020inverse}.  In comparison to other methods, such as
finite-difference time-domain (FDTD) methods, semi-analytical methods often
have much lower computational cost.  

This advantage stems from how semi-analytical methods discretize
Maxwell's equations.  In contrast to other approaches (e.g., FDTD methods) that discretize 
the spatial domain fully (i.e., in
all three dimensions)~\cite{yee1966numerical}, semi-analytical methods discretize the
spatial domain partially (e.g., in only $x$- and $y$-dimension but not
$z$-dimension).  This is possible because many photonic structures have a
primary light propagation direction (referred in this paper as $z$-direction;
see \figref{notation_diagram}). In some cases, along the light propagation direction, the
structure's cross-sectional shape stays unchanged (e.g., a metasurface unit).
Therefore, we do not have to discretize the structure along $z$-direction;
instead, light propagation in the structure can be viewed as superposition of
individual propagating modes experiencing phase shifts.  This is the
fundamental view that enables semi-analytical methods to reduce computational
cost (see \textbf{2.1}).

\begin{figure}[t]
    \centering
    \includegraphics[width=0.95\textwidth]{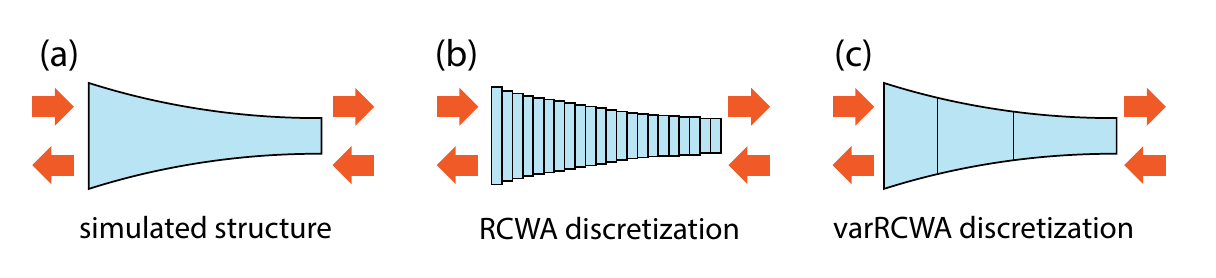}
    \vspace{-5mm}
    \caption{\textbf{Discretization in conventional and our method.} 
    Provided a photonic structure with varying cross-sectional shape (a), 
    conventional semi-analytical methods (such as RCWA) must discretize the structure along the wave propagation direction
    into many thin sections (b), leading to expensive computation.
    (c) Our method discretizes the structure adaptively into a much smaller number of sections, thanks 
    to the use of higher-order perturbative expansions. 
    The arrows indicate the incoming and outgoing direction of light waves.}
    \label{fig:discretization}
\end{figure}

However, this view becomes unsound for many photonic structures wherein along
the primary light propagation direction, the structure's cross section
varies~\cite{piggott2015inverse, miller2020large}.  A common example is
photonic waveguides (such as a taper; see \figref{discretization}-a). To
simulate these photonic structures using semi-analytical methods, one has to
further discretize the structure along $z$-direction into a series of thin
sections~\cite{liu2012s4,jing2013analysis} (see \figref{discretization}-b).
In each section, the cross-sectional shape is assumed unchanged, and thereby a
semi-analytical method can be used to simulate that section. Yet, this approach
requires a large number of sections, which in turn devastate the computational
advantage of semi-analytical methods.  Apart from the computational cost, it is
often unclear how many discrete sections are needed to achieve certain accuracy.
In practice, one has to rely on trail and error to choose a proper resolution
for sufficient accuracy.  Oftentimes, to obtain satisfactory results, multiple
runs of the simulation method (each with a different resolution) are needed.


In this work, we overcome these limitations.  
Our method requires no trial and error, thus much easier to use: provided a photonic structure 
and a user-specified accuracy level (i.e., a real number), our method
automatically decides how to discretize the structure in $z$-direction, aiming
to reduce the overall computational cost while achieving the desired accuracy.
To obtain simulation results of user-specified accuracy, only one run is needed.

To this end, our core development is twofold:
1) We show that the conventional semi-analytical methods (such as RCWA) are merely zeroth-order
approximation with respect to the structure's cross-sectional variation. 
Through a novel change of variable,
we propose a high-order semi-analytical method, which 
allows the structure's cross section to vary over $z$-direction,
without discretizing it into thin sections.
2) Leveraging this high-order method, we introduce an algorithm 
that automatically and adaptively discretizes the 
structure to achieve a user specified accuracy level.
For regions where the
cross section varies rapidly in $z$-direction, our algorithm will slice the structure in fine resolution 
to ensure simulation accuracy; for regions with little
cross-sectional variance, it will discretize them coarsely to save computational
cost. 

%
%

We use our method to analyze various photonic structures, and compare it with  
conventional semi-analytical methods (such as RCWA).
We show that our method, as a higher-order approach, indeed converges faster.
As a result, to obtain the same level of accuracy, our method requires much less
computational time and no resolution tuning at all.

\section{Method}
%
We now present our core development.  To understand the rationale behind our
development, we start by briefly reviewing the limitations of widely used
semi-analytical methods.


\subsection{Established Semi-analytical Methods and Their Limitations}  
Semi-analytical methods discretize the spatial domain of Maxwell's equations in $x$- and $y$-directions
but not in $z$-direction. In frequency domain, Maxwell's equations become into
\begin{equation}\label{eq:pq_form}
    \frac{\partial \bvec{e}}{\partial z} = \frac{j}{k_0}\mat{P}\bvec{h}
    \;\textrm{ and }\;
    \frac{\partial \bvec{h}}{\partial z} = \frac{j}{k_0}\mat{Q}\bvec{e},
\end{equation}
where the vectors $\bvec{e}$ and $\bvec{h}$ are discrete representations of the electric 
and magnetic fields on an $xy$-plane at a $z$ position, and the matrices $\mat{P}$ and $\mat{Q}$
encode the distributions of material permeability and permittivity on the $xy$-plane at the same $z$ position.
Depending on specific representations of $\bvec{e}$ and $\bvec{h}$, different semi-analytical methods emerge.
The most widely used (e.g., for the analysis of 
metasurfaces~\cite{divitt2019ultrafast}, gratings~\cite{mohamad2020fast} and waveguides~\cite{zhu2020inverse}) is  
rigorous coupled wave analysis (RCWA) method, wherein $\bvec{e}$ and $\bvec{h}$ are discretized using 2D Fourier
basis on the $xy$-plane. In this paper, 
our development can be applied to semi-analytical methods in general (such as the method of lines~\cite{pregla1989method}),
although our implementation and numerical experiments focus in particular on the RCWA method.

\begin{figure}[t]
    \centering
    \includegraphics[width=0.89\textwidth]{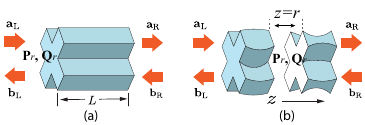}
    \vspace{-1mm}
    \caption{\textbf{Notation.} Provided a photonic structure section,
    we model light wave propagation in four components: the incoming and outgoing 
    waves at both the left and right ends, as shown by the orange arrows.
    (a) A section with constant cross-sectional shape has 
    fixed material matrices $\mat{P}_r$ and $\mat{Q}_r$. 
    (b) When a section has a varying cross-sectional shape, we choose a specific position $z=r$ (called reference 
    position), and use the material matrices, $\mat{P}_r$ and $\mat{Q}_r$, at $z=r$ to construct the basis 
    for describing forward and backward waves. Our method treats such a section as a perturbation of a section
    with constant cross-sectional shape in (a), and leverages a perturbative expansion 
    for numerical simulation. Here for visualization purpose, we cut through the section at $z=r$ to reveal
    the cross-section at $z=r$, which may differ from the cross-sectional shape at its two ends.}
    \label{fig:notation_diagram}
\end{figure}

When the structure's cross-sectional shape is fixed along $z$-direction, both
$\mat{P}$ and $\mat{Q}$ in \eq{pq_form} are constant matrices, and
\eqnref{pq_form} can be solved through an eigenvalue decomposition, that is, $\mat{P}\mat{Q}=\mat{W}\mat{\Lambda}^{2}\mat{W}^{-1}$.
The resulting $\mat{W}$ and $\mat{\Lambda}$
allow us to express the solution of \eq{pq_form} as
\begin{equation}\label{eq:eh_sol}
    \bvec{e} = \frac{1}{2}\mat{W}\left(e^{\mat{\Lambda}z}\bvec{a}_L + e^{-\mat{\Lambda}z}\bvec{b}_L\right)
    \quad \text{and} \quad
    \bvec{h} = \frac{1}{2}\mat{V}\left(e^{\mat{\Lambda}z}\bvec{a}_L - e^{-\mat{\Lambda}z}\bvec{b}_L\right).
\end{equation}
Here $\mat{V}$ is the basis for describing the cross-sectional magnetic field, related to the
eigenvectors $\mat{W}$ (i.e., the basis for electric field) through $\mat{V}=\mat{Q}\mat{W}\mat{\Lambda}^{-1}$;
$\bvec{a}_{\textrm{L}}$ and $\bvec{b}_{\text{L}}$ are vectors stacking the coefficients of the forward and backward light waves
at the {left end} of the structure (see \figref{notation_diagram}). In
addition, from \eqnref{eh_sol}, we can define the structure's \emph{scattering
matrix}, which relates the output state of an optical wave after propagating
through the structure with its input state, namely
\begin{equation}
\begin{bmatrix}
 \bvec{a}_{\text{R}} \\ 
 \bvec{b}_{\text{L}}
\end{bmatrix} = \mat{S}
\begin{bmatrix}
 \bvec{a}_{\text{L}} \\ 
 \bvec{b}_{\text{R}}
\end{bmatrix},\quad\text{where}\;\;
\mat{S} = \begin{bmatrix}
    \mat{T}_{\text{LR}} & \mat{R}_{\text{R}} \\
    \mat{R}_{\text{L}} & \mat{T}_{\text{RL}} 
\end{bmatrix},
\end{equation}
where, corresponding to $\bvec{a}_{\text{L}}$ and $\bvec{b}_{\text{L}}$,
$\bvec{a}_{\text{R}}$  and $\bvec{b}_{\text{R}}$ describe the forward and
backward waves at the {right end} (see \figref{notation_diagram}).  Once the scattering matrix is known,
the structure's optical performance (e.g., mode conversion efficiency and phase
shift of a waveguide) can be directly computed.

When the structure's cross-sectional shape varies along $z$-direction, one has to
split the structure into a series of small sections so that every section can be
approximated as having a fixed cross section (\figref{discretization}-b).  Each
section $i$ is then analyzed through the aforementioned process, from which one
can compute its scattering matrix $\mat{S}_i$. Finally, by
combining all the scattering matrices using Redheffer star
product~\cite{redheffer1959inequalities}, the entire structure's scattering matrix is obtained.

\paragraph{Limitations.}
In semi-analytical methods, the assumption of having a fixed cross section can be viewed as a zeroth-order
approximation of the structure (as derived in \secref{high_order}). 
To achieve sufficient accuracy, such a crude approximation must be remedied with {small structure length}.
As a result, even mild cross-sectional variation requires a large number of
discrete sections (\figref{discretization}-b).
For every section, an eigenvalue decomposition is needed,
and thus its computational cost is further scaled by the total number of sections.
To reduce the computational cost, we need to
reduce the total number of discrete sections (\figref{discretization}-c) while retaining simulation accuracy.
This motivates us to seek a high-order semi-analytical method, one that accounts for
the cross-sectional variation in a long section and thereby reduces the total number of sections.

\subsection{High-order Semi-analytical Methods}\label{sec:high_order}
Consider a section of photonic structure along $z$-direction. Suppose its cross-sectional shape varies, that is, 
in \eqnref{pq_form}, $\mat{P}$ and $\mat{Q}$ are not constant matrices; they change over $z$.
In this case, the solution of \eq{pq_form} is not as simple as \eq{eh_sol}.
Now, our goal is to express the solution of  \eq{pq_form} as a {perturbative}
expansion with respect to cross-section variation (\figref{notation_diagram}), and this expansion will serve as the core numerical recipe
of our method.

\paragraph{Na\"{i}ve solution.} To understand the insight of our development, we start with a na\"{i}ve (but impractical) expansion 
form of the solution. 
First, inspired by \eqnref{eh_sol}, we use a set of basis vectors $\mat{W}$ and $\mat{V}$ to describe cross-sectional
electric and magnetic fields respectively\textemdash the specific choice of $\mat{W}$ and $\mat{V}$
in presence of varying cross section will be described shortly. The cross-sectional electric and magnetic fields,
$\bvec{e}$ and $\bvec{h}$, are formed by light waves propagating forward and backward in the structure, 
with the relations:
\begin{equation}\label{eq:forward_backward}
    \bvec{a}(z) = \mat{W}^{-1}\bvec{e}(z) + \mat{V}^{-1}\bvec{h}(z)
    \;\textrm{ and }\;
    \bvec{b}(z) = \mat{W}^{-1}\bvec{e}(z) - \mat{V}^{-1}\bvec{h}(z),
\end{equation}
where $\bvec{a}$ and $\bvec{b}$ are coefficients in the chosen basis for describing the forward and backward waves, respectively,
and they vary over $z$. Next, to establish a differential equation of $\bvec{a}$ (and $\bvec{b}$), we differentiate both sides
of~\eq{forward_backward}, and then using Maxwell's equations~\eq{pq_form}, we obtain
\begin{equation}\label{eq:direct_u_deriv}
\frac{\partial \bvec{a}}{\partial z} =\mat{W}^{-1}\frac{\partial \bvec{e}}{\partial z} + \mat{V}^{-1}\frac{\partial \bvec{h}}{\partial z}
= \frac{j}{2k_0}\left[\mat{W}^{-1}\mat{P}\mat{V}\left(\bvec{a} - \bvec{b}\right) + \mat{V}^{-1}\mat{Q}\mat{W}\left(\bvec{a} + \bvec{b}\right)\right].
\end{equation}
Equivalently, we have the integral equation
\begin{equation}\label{eq:int_eqn}
    \bvec{a}(z) = \bvec{a}(z_{L}) + \frac{j}{2k_0}\int^{z}_{z_{L}}\left[\left(\mat{W}^{-1}\mat{P}\mat{V}+\mat{V}^{-1}\mat{Q}\mat{W}\right)\bvec{a}(z')+\left( \mat{V}^{-1}\mat{Q}\mat{W}-\mat{W}^{-1}\mat{P}\mat{V}\right)\bvec{b}(z')\right]\diff z',
\end{equation}
where $z_L$ is the starting $z$ position of the considered structure section.
Note also that both $\mat{P}$ and $\mat{Q}$ vary over $z$, and a similar
integral equation can be obtained for $\bvec{b}(z)$. Equation~\eqref{eq:int_eqn}, in
theory, allows us to express $\bvec{a}(z)$ as a {perturbative} expansion. 
This is achieved by recursively substituting $\bvec{a}(z')$ in the integrand with \eqnref{int_eqn}
itself up to a certain order (and similarly for $\bvec{b}(z)$).
For example, to obtain a first-order expansion, one can replace $\bvec{a}(z')$ and $\bvec{b}(z')$ in \eq{int_eqn} with their 
zeroth-order approximations $\bvec{a}_L \approx \bvec{a}(z_L)$ and $\bvec{b}_R \approx \bvec{b}(z_R)$.


To use this expansion for analyzing a long section (and thereby reduce the total number
of sections), the norms of $\mat{W}^{-1}\mat{P}\mat{V}$ and
$\mat{V}^{-1}\mat{Q}\mat{W}$ must be sufficiently small\textemdash an intuitive
explanation of this requirement is provided in Supplement 1.  This
requirement, however, is hardly satisfied in practice, as both $\mat{P}$ and
$\mat{Q}$ depend on the cross-sectional material distribution, and their norms may
become arbitrarily large.  Nevertheless, the development of this expansion
motivates a viable strategy: 
in order to obtain a stable perturbative expansion, 
we need to avoid using $\mat{P}$ and $\mat{Q}$ in an integral equation like~\eq{int_eqn}; instead, we seek an expansion that
involves only the variation of $\mat{P}$ and $\mat{Q}$ over $z$.

\paragraph{Preconditioned solution.}
Our proposed expansion starts with a change of variables. We introduce two variables:
\begin{equation}\label{eq:u_d_ph}
    \tilde{\bvec{a}} = e^{-\frac{j}{k_0}\mat{\Lambda}z}\bvec{a}
    \;\textrm{ and }\;
    \tilde{\bvec{b}} = e^{\frac{j}{k_0}\mat{\Lambda}z}\bvec{b},
\end{equation}
where $\mat{\Lambda}$ is the eigenvalue matrix resulted from eigen-decomposition
$\mat{P}_r \mat{Q}_r=\mat{W}\mat{\Lambda}^2\mat{W}^{-1}$.
Here, $\mat{P}_r$ and $\mat{Q}_r$ are fixed matrices encoding the distribution of material
permeability and permittivity at a particular $z=r$ position. 
Ideally, the cross section at $z=r$ is chosen to represent the 
``average'' cross section over the entire section so that
it can be used to construct the basis vectors $\mat{W}$ and $\mat{V}$.
We therefore refer to this position as the \emph{reference position}
of the section (see \figref{notation_diagram}-b). 
While in theory one can choose any $r$ position, in practice we simply use the
mid-point of the section.  The resulting $\mat{W}$ (and
$\mat{V}$ through $\mat{V}=\mat{QW\Lambda}^{-1}$) is used as the basis for describing
forward and backward waves (recall~\eqnref{forward_backward}).

This change of variables is the key to introduce the variation 
of $\mat{P}$ and $\mat{Q}$ in an integral equation similar to~\eq{int_eqn}.
By differentiating~\eq{u_d_ph} and 
using \eqnref{pq_form}, we obtain the following differential equations (see the derivation in Supplement 2):
\begin{equation} \label{eq:diff_ab} 
    \begin{split}
    \frac{\partial \tilde{\bvec{a}}}{\partial z} &= \frac{j}{2k_0}e^{\left(-\frac{j}{k_0}\mat{\Lambda}z\right)}\delta \mat{A}e^{\left(\frac{j}{k_0}\mat{\Lambda}z\right)} \tilde{\bvec{a}} - \frac{j}{2k_0}e^{\left(-\frac{j}{k_0}\mat{\Lambda}z\right)}\delta \mat{B}e^{\left(-\frac{j}{k_0}\mat{\Lambda}z\right)} \tilde{\bvec{b}},\\
    \frac{\partial \tilde{\bvec{b}}}{\partial z} &= \frac{j}{2k_0}e^{\left(\frac{j}{k_0}\mat{\Lambda}z\right)}\delta \mat{B}e^{\left(\frac{j}{k_0}\mat{\Lambda}z\right)} \tilde{\bvec{a}} - \frac{j}{2k_0}e^{\left(\frac{j}{k_0}\mat{\Lambda}z\right)}\delta \mat{A}e^{\left(-\frac{j}{k_0}\mat{\Lambda}z\right)} \tilde{\bvec{b}},
    \end{split}
\end{equation}
where $\delta \mat{A}$ and $\delta \mat{B}$ are related to the material variation along $z$-direction, namely, 
\begin{equation} \label{eq:delta_AB} 
    \begin{split}
    \delta \mat{A}(z) = \mat{W}^{-1}\left(\mat{P}(z)-\mat{P}_r\right)\mat{V}+\mat{V}^{-1}\left(\mat{Q}(z)-\mat{Q}_r\right)\mat{W}, \\
    \delta \mat{B}(z) = \mat{W}^{-1}\left(\mat{P}(z)-\mat{P}_r\right)\mat{V}-\mat{V}^{-1}\left(\mat{Q}(z)-\mat{Q}_r\right)\mat{W}.
    \end{split}
\end{equation}

Next, we rewrite \eqnref{diff_ab} in integral forms and replace $\tilde{\bvec{a}}$ and $\tilde{\bvec{b}}$ 
using \eqnref{u_d_ph}. This leads to a new set of integral equations of $\bvec{a}(z)$ and $\bvec{b}(z)$,
ones that differ from \eq{int_eqn}:
\begin{align}
  \bvec{a}(z) &= e^{\frac{j}{k_0}\mat{\Lambda}(z-z_L)}\bvec{a}_L+\frac{j}{2k_0}\int^{z}_{z_L}e^{\frac{j}{k_0}\mat{\Lambda}(z-z')}\delta\mat{A}(z')\bvec{a}(z')\diff z' - \frac{j}{2k_0}\int^{z}_{z_L}e^{\frac{j}{k_0}\mat{\Lambda}(z-z')}\delta \mat{B}(z')\bvec{b}(z')\diff z' \label{eq:u_integral} \\
  \bvec{b}(z) &= e^{\frac{j}{k_0}\mat{\Lambda}(z_R-z)}\bvec{b}_R-\frac{j}{2k_0}\int^{z_R}_{z}e^{\frac{j}{k_0}\mat{\Lambda}(z' - z)}\delta\mat{B}(z')\bvec{a}(z')\diff z' + \frac{j}{2k_0}\int^{z_R}_{z}e^{\frac{j}{k_0}\mat{\Lambda}(z' - z)}\delta \mat{A}(z')\bvec{b}(z')\diff z'. \label{eq:d_integral}
\end{align}
Note that here the two equations are integrated from different directions. In \eqnref{u_integral},
we use $\bvec{a}_L$ at the left end of the section as the initial value,
and the integral is along the forward direction. On the contrary, \eqnref{d_integral}
integrates along the backward direction, using $\bvec{b}_R$ as the initial value.
In this way, all the propagation phase terms (i.e., $e^{j\mat{\Lambda}\Delta
z/k_0}$) become smaller than one, and thereby the numerical computation of these
integrals stays stable.

Now, we can construct {perturbative} expansions of the outgoing waves $\bvec{a}_R$
and $\bvec{b}_L$.  Similar to the construction of the na\"{i}ve expansions
above, we substitute $\bvec{a}(z')$ and $\bvec{b}(z')$ in the integral terms
with Eqs.~\eq{u_integral} and~\eq{d_integral} themselves, and this substitution
is done recursively up to a certain order. These expansions offer a numerical
recipe for analyzing how the structure interacts with propagating waves.
For example, provided the input waves to the section, including the forward wave at the left end 
(described by $\bvec{a}_L$; see \figref{notation_diagram})
and the backward wave at the right end (described by $\bvec{b}_R$), the outgoing waves at the right and left ends 
can be estimated using the following first-order expansions:
\begin{dmath}\label{eq:u_pert_1}
    \bvec{a}_R \approx e^{\frac{j}{k_0}\mat{\Lambda}(z_R-z_L)}\bvec{a}_L+\frac{j}{2k_0}\int^{z_R}_{z_L}e^{\frac{j}{k_0}\mat{\Lambda}(z_R-z')}\delta\mat{A}(z')e^{\frac{j}{k_0}\mat{\Lambda}(z'-z_L)}\diff z'\bvec{a}_{L} - \frac{j}{2k_0}\int^{z_R}_{z_{L}}e^{\frac{j}{k_0}\mat{\Lambda}(z_R-z')}\delta \mat{B}(z')e^{\frac{j}{k_0}\mat{\Lambda}(z_R-z')}\diff z'\bvec{b}_R,
\end{dmath}
\vspace{-4mm}
\begin{dmath}\label{eq:d_pert_1}
    \bvec{b}_L\approx e^{\frac{j}{k_0}\mat{\Lambda}(z_R-z_L)}\bvec{b}_R+ \frac{j}{2k_0}\int^{z_R}_{z_L}e^{\frac{j}{k_0}\mat{\Lambda}(z'-z_L)}\delta \mat{A}(z')e^{\frac{j}{k_0}\mat{\Lambda}(z_R-z')}\diff z'\bvec{b}_R-\frac{j}{2k_0}\int^{z_R}_{z_L}e^{\frac{j}{k_0}\mat{\Lambda}(z'-z_L)}\delta\mat{B}(z')e^{\frac{j}{k_0}\mat{\Lambda}(z'-z_L)}\diff z'\bvec{a}_L.
\end{dmath}

Remarkably, these expansions do not involve $\mat{W}^{-1}\mat{P}\mat{V}$
or $\mat{V}^{-1}\mat{Q}\mat{W}$ (unlike \eq{int_eqn}).
But rather, their integrands depend on $\delta\mat{A}$ and  $\delta\mat{B}$.
Thus, Eqs.~\eq{u_pert_1} and~\eq{d_pert_1} can be viewed as a (first-order) perturbation 
solution of the outgoing waves.
If the cross section is invariant, 
$\delta\mat{A}=\delta\mat{B} = \mat{0}$, and the integral terms 
in Eqs.~\eq{u_pert_1} and~\eq{d_pert_1} vanish.
In this case, the expressions are precisely the same as the conventional 
semi-analytical methods, indicating that the conventional semi-analytical
methods are zeroth-order perturbation in presence of cross-sectional variation.
If the cross section is slowly varying over $z$, the norms of $\delta\mat{A}$ and  $\delta\mat{B}$
are small, and thereby the first-order perturbation converges even for a long section.
This contrasts starkly to the na\"{i}ve expansions based on \eqnref{int_eqn}, whose 
reliance on $\mat{P}$ and $\mat{Q}$ drastically restricts the section length.
If the cross section varies quickly,
the norms of $\delta\mat{A}$ and  $\delta\mat{B}$ may be large,
and we have to use short section length to ensure convergence of our 
perturbation solution.

\textbf{Discussion: alternative approaches.}
Maxwell's equation~\eqref{eq:pq_form} may also be viewed as an initial value problem provided
with the electric and magnetic fields at $z=z_L$.
Solution to this initial value problem can be expressed as a product integral
that involves matrix exponentials~\cite{helton1976numerical, slavik2007product}:
\begin{equation}\label{eq:product_integral_form}
    \begin{bmatrix} \bvec{e}(z_R) \\ \bvec{h}(z_R)\end{bmatrix}= \lim_{N\to\infty}\prod^{1}_{i=N} \exp{\left(\frac{j(z_i - z_{i-1})}{k_0}\begin{bmatrix} \mat{0}& \mat{P}_r(z_{i-1}\rightarrow z_i)\\
        \mat{Q}(z_{i-1}\rightarrow z_i)& \mat{0}\end{bmatrix}\right)}\begin{bmatrix}\bvec{e}(z_L) \\ \bvec{h}(z_L)\end{bmatrix}.
\end{equation}
In this expression, the structure is also discretized into a series of sections at $z_i$,
where $z_L = z_0 < z_1 < \cdots < z_N = z_R$; 
and $\mat{P}_r$ and $\mat{Q}_r$ are permeability and permittivity matrices
evaluated at the reference position of each section. 

Numerical evaluation of~\eqref{eq:product_integral_form}, however, is rather challenging.
Prior works use a low-order Taylor expansion of the matrix exponentials to
evaluate~\eqref{eq:product_integral_form}~\cite{roberts2017rigorous, li2020efficient}.
But to use this expansion, section length must be excessively short (i.e., $z_i -z_{i-1}\le0.1\lambda$
where $\lambda$ is the wavelength), and a large number of sections are needed.

Another approach is to convert the product of matrix exponential into the
exponential of matrix summation, similar to $\prod\exp{\left(x_i\right)} =
\exp{\left(\sum_i x_i\right)}$ for scalar values $x_i$.  However, this is not
straightforward because matrix multiplication is noncommutative. As a result, 
Magnus expansion ~\cite{blanes1998magnus, iserles1999solution,
blanes2009magnus} and Fer expansion~\cite{takegoshi2015comparison} have been
used to correct the use of exponential of matrix summation. In this vein, RCWA
can be viewed as an approximation of Magnus expansion, as shown
in~\cite{chu2003finite}.  Although a high-order approximation of Magnus
expansion was introduced in~\cite{chu2003finite}, it still requires each
discrete section to be short even with small cross-sectional variation.

A more fundamental problem of this line of approaches stems from the treatment of
initial value problem.  Light propagating in photonic structures 
almost always scatters, resulting in forward and backward propagations. But treating it as an initial value problem 
implies that only a single direction is considered, thus causing the numerical integration unstable.  It is this
reason that our method, like conventional RCWA, separates the light waves into
forward and backward going components, and integrate them along two separate
directions (recall \eqref{eq:u_integral} and \eqref{eq:d_integral}).

\subsection{Numerical Implementation with Adaptive Discretization}
\paragraph{Numerical integration.}
The integrals in Eqs.~\eq{u_pert_1} and~\eq{d_pert_1} can be numerically evaluated using a quadrature rule.
In practice, we use the quadrature rule that samples three positions, at $\zL$, $\frac{\zL + \zR}{2}$, and $\zR$, respectively,
and the integral is estimated by weighted summation of the integrand values at the sampled positions.
This quadrature rule involves only matrix multiplications, and thus can be
computed at low cost; in practice, we implement it on a Graphics Processing Unit (GPU). 
Consequently, the major cost of evaluating wave propagation in a single section (using~\eq{u_pert_1}
and~\eq{d_pert_1}) comes from the eigen-decomposition that computes $\mat{W}$, $\mat{\Lambda}$, and $\mat{V}$
at the reference position. For a single structure section, 
the eigen-decomposition is also needed in conventional semi-analytical
methods. Our perturbative expansion is more accurate at a trivial cost of the
additional numerical integration. But when it comes to simulating an
entire photonic structure, our expansion enables adaptively discretization of the
structure along $z$-direction, thereby reducing the overall computational
cost and outperforming the conventional methods.
Next, we describe our algorithm for simulating an entire structure. 

\begin{algorithm}[t]
\DontPrintSemicolon
    \KwInput{Simulated region along z axis: $[z_\min, z_\max]$, error bound $\alpha$}
    \KwOutput{Scattering matrix from $z_\min$ to $z_\max$}
    \KwData{Simulated Geometry $\mathcal{G}$}
    Evaluate high order scattering matrix $S$ according to Supplement 3 by sampling on $\mathcal{G}$\\
    \If{$\text{estimated error}$ (see Supplement 6)~$ < \alpha$ }
    {
        \Return $S$
    }
    \Else
    {
        Subdivide this section into $M$ subsections $(z_0, z_1)$, $(z_1, z_2)$, $\dots$,  $(z_{M-1}, z_M)$ \\
        ~where $z_\min = z_0 < z_1 < \dots < z_M = z_\max$\\
        Find a reference point for each subsection. \\
        \For {$i = 1$, ..., $M$}
        {
            scattering matrix $S_i(z_{i-1}\to z_i) \gets$  Adaptive-Variant-RCWA($z_{i-1}$, $z_i$, $\alpha$)\\
            Project the basis of $S_i$ to match the basis of $S_{i-1}$ (see Supplement 4)
        }
        \Return $S \leftarrow $ Redheffer star product of $S_1$, ...,  $S_M$.
    }
\caption{Adaptive-Variant-RCWA($z_\min$, $z_\max$, $\alpha$)}\label{alg:adaptive_variant_rcwa}
\end{algorithm}

\paragraph{Adaptive discretization.}
%
We now consider the simulation of an entire photonic structure, not just a
single section.  We propose an algorithm that adaptively discretizes the
structure in $z$-direction and simulates wave propagation.  Input to our
algorithm is a photonic structure and a desired accuracy level (i.e., a real
number).  The goal here is to achieve the desired accuracy while reducing the
overall computational cost.

Our proposed algorithm is outlined in \algref{adaptive_variant_rcwa}.
The key idea behind this algorithm is as follows:
starting from the entire structure as a single section, it recursively
subdivides a structure section into $M$ subsections (Line 5 in
\algref{adaptive_variant_rcwa}).  The subdivision occurs when an estimated
simulation error of the current section is larger than the user specified
accuracy level $\alpha$ (Line 2 in \algref{adaptive_variant_rcwa}).
Afterwards, all subsections are simulated individually, and they may be subdivided
further in a recursive way (Line 9).  For each subsection, we compute how
the wave propagates using~\eq{u_pert_1} and~\eq{d_pert_1}, and further compute this subsection's 
scattering matrix (Line 1; see Supplement 3 for detailed formulas).
Lastly, the scattering matrices of individual subsections are combined through Redheffer star
product~\cite{redheffer1959inequalities} to form the scattering matrix of the
parent section (Line 11).


\begin{figure}[t]
    \centering
    \vspace{-10mm}
    \includegraphics[width=0.9\textwidth]{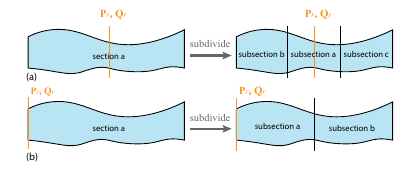}
    \vspace{-5mm}
    \caption{(a) we choose the midpoint of each section as 
        its reference point to evaluate $\mat{P}_r$ and $\mat{Q}_r$.
When subdividing one section, we subdivide it evenly into three subsections.
In this way, the midpoint of the mid-section is the same as that of the parent section, and thus
$\mat{P}_r$ and $\mat{Q}_r$ can be reused.
(b) We also experimented with subdividing a section evenly into two subsections (i.e., $M=2$). 
This strategy is proper when we use the endpoint of a section as its reference point.}
    \label{fig:two_ways_of_subdivision}
\end{figure}

We estimate the error by computing the discrepancy of asymptotic expansions in two
consecutive orders (e.g., between zeroth order and first order) Supplement 6\textemdash
such an error estimation has been used in other asymptotic expansions~\cite{Chen14:ANM,cochelin1994path}.
In practice, when a section is subdivided, we subdivide it into three
subsections (i.e., $M=3$).  In this way, the eigen-decomposition needed for the
parent section can be reused for simulating the central subsection (see
\figref{two_ways_of_subdivision}-b). An alternative is to use binary subdivision (i.e., $M=2$), which
we will compare in our numerical experiments in \textbf{3.2}.

Through recursive subdivision, adaptive discretization on $z$-direction
naturally emerges: regions with rapidly varying cross-sections will be more
subdivided, and thereby the algorithm automatically uses fine sections to
ensure accuracy; meanwhile, in smoothly varying regions, it uses low
resolution to save computation.

\section{Results}
We conduct numerical experiments to validate the accuracy of our method and
compare its performance with conventional semi-analytical methods. Since there
are different ways of implementing semi-analytical methods (depending on how
$\bvec{e}$ and $\bvec{h}$ in \eq{pq_form} are represented), we adopt the most
widely used representation, the RCWA method~\cite{moharam1981rigorous}.

We implement both our method and the conventional RCWA in C++ programming language.
Numerical computation in both methods, such as matrix multiplication and eigen-decomposition,
can benefit from parallel computation. We therefore leverage CUDA~\cite{cuda} on Nvidia Graphics Processing Units (GPUs)
to accelerate the computation in both methods.

\begin{figure}[!htp]
    \centering
    \includegraphics[width=0.99\textwidth]{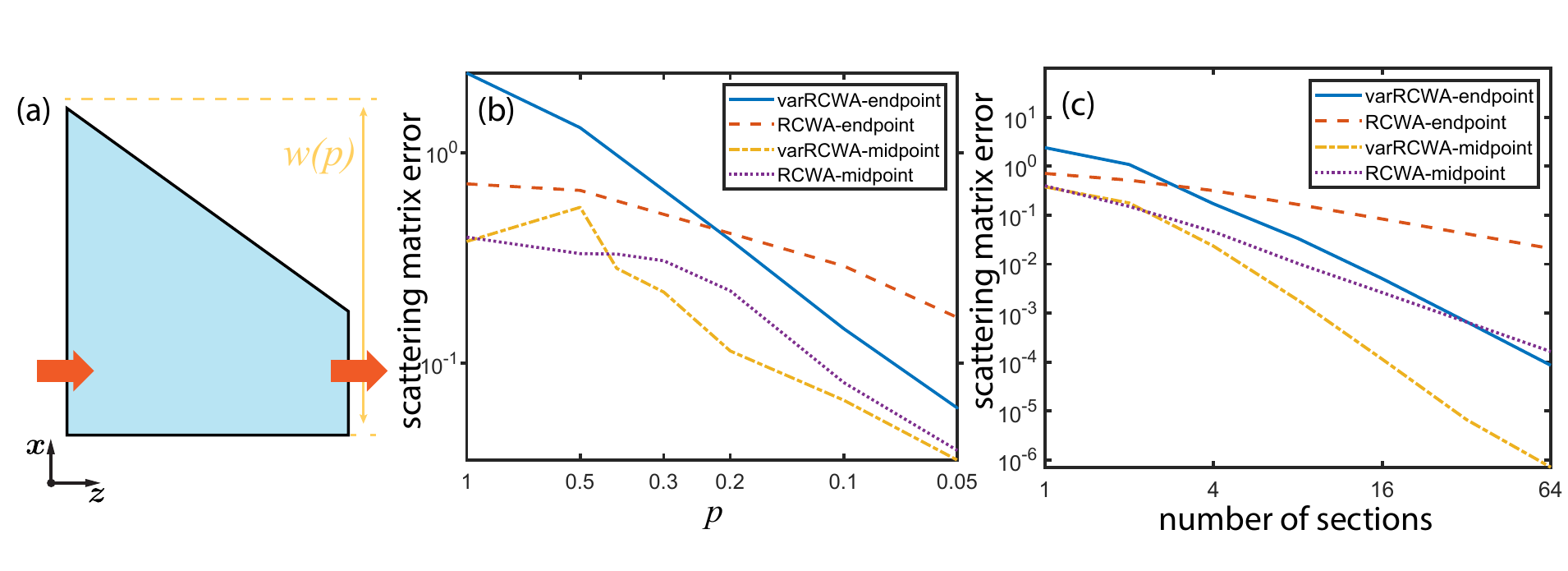}
    \vspace{-4.5mm}
    \caption{
(a) We use VarRCWA to analyze a waveguide whose cross-sectional shape changes linearly.
The length of the waveguide is $1 ~\um$. 
Its width on the right end is fixed at $2.6~\um$, and we use a parameter $p$ to specify 
the width on the left end. When $p=1$, the left width is $3.7~\um$, and when $p=0$, the 
left width is the same as the right end (i.e., $2.6~\um$).
(b) We simulate this waveguide with no $z$-direction discretization (i.e., treating 
it as a single section), and measure the scattering matrix errors as $p$ changes.
(c) We then increase the $z$-direction resolution,
and measure the scattering matrix errors as $p$ changes.
In the legends, 
``midpoint'' indicates the use of the middle position as the reference position
to construct $\mat{P}_r$ and $\mat{Q}_r$(see \figref{two_ways_of_subdivision}-a),
and ``endpoint'' indicates the use of the end position of a section as the reference position (see \figref{two_ways_of_subdivision}-b). See \textbf{Data File 1-2} for underlying values.
}\label{fig:high_order_validation}
\end{figure}
\subsection{Validation}
To validate the accuracy of our method, we consider a trapezoid-shaped
waveguide as shown in \figref{high_order_validation}-a. This waveguide has a
fixed $220 {\rm nm}$ thickness, and we use a parameter $p$ to control its cross-sectional
variation from its left to right end (see \figref{high_order_validation} caption).

First, we verify that our method has higher accuracy than the conventional method. 
To evaluate the accuracy, we first discretize the waveguide into 256 sections,
and compute the entire waveguide's scattering matrix $\mat{S}^*$ using the
conventional RCWA. This $\mat{S}^*$ is then used as a ground-truth to measure a
max norm error, $\|\mat{S}-\mat{S}^*\|_\max$, where $\mat{S}$ is the scattering
matrix resulted from either our method or the conventional RCWA.  In both
cases, we use a single section to simulate the entire waveguide for fair
comparison. We measure the errors while the cross-sectional variation (controlled by $p$)
changes. As shown in \figref{high_order_validation}-b, when the cross-sectional variation 
is large (i.e., $p\to 1$), our method is more accurate than conventional RCWA, thanks to its higher-order approximation.
When the variation becomes small (i.e., $p\to 0$), 
both methods have comparable accuracy.



Further, we evaluate the convergence of our method. 
We again consider the waveguide shown in \figref{high_order_validation} (with $p=1$),
and evaluate the scattering matrix errors of our method and the conventional RCWA
while the $z$-direction resolutions (in both methods) increases.
The results are reported in \figref{high_order_validation}-c: our method converges faster;
as the resolution becomes higher, our method is much more accurate. 


\begin{figure}[!htp]
    \vspace{-9mm}
    \centering
    \includegraphics[width=0.95\textwidth]{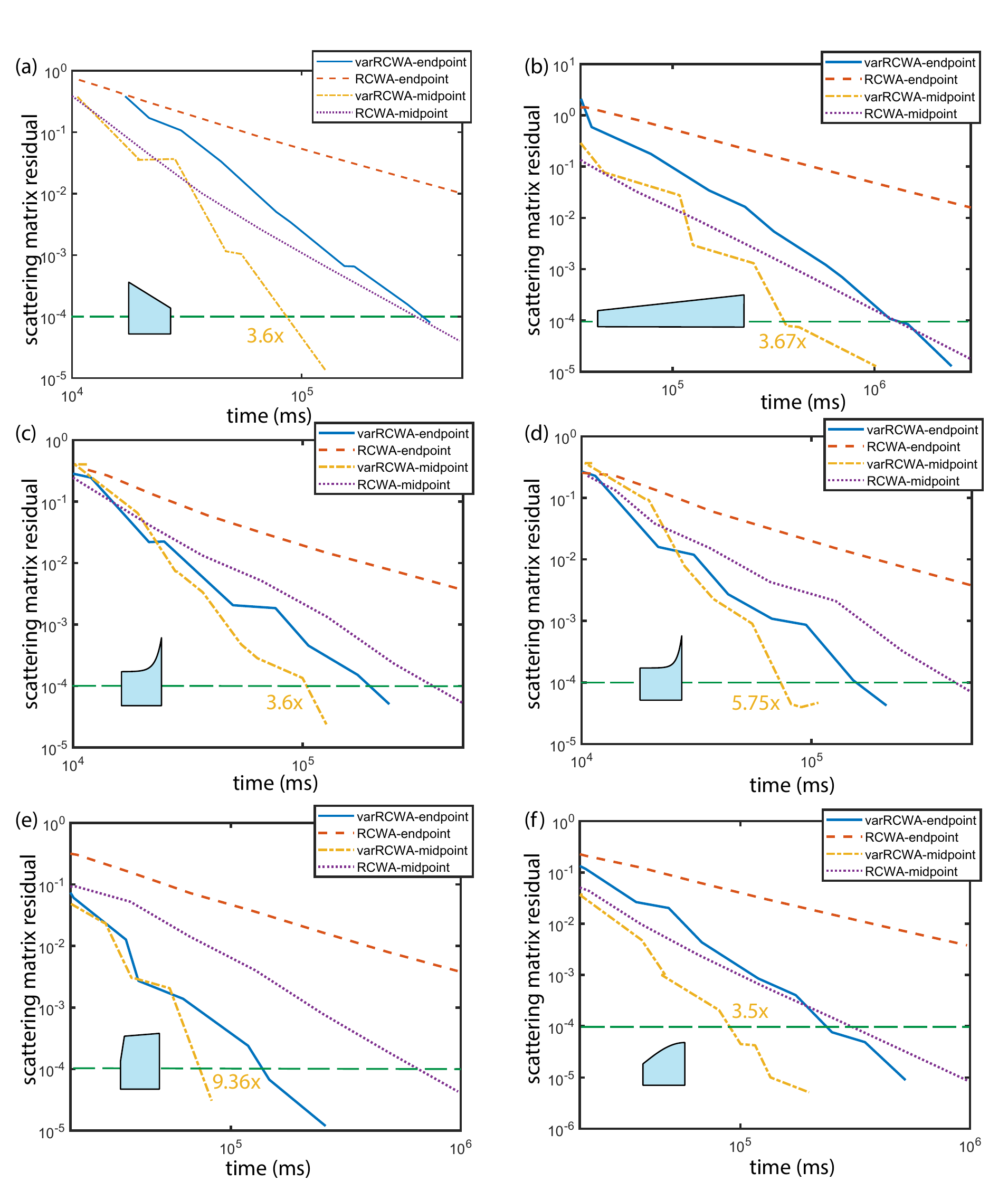}
    \vspace{-4mm}
    \caption{\textbf{Tests on various waveguide geometries.} 
(a) a trapezoid waveguide. (b) a
long waveguide of $10 \um$ length. (c-d) waveguides with exponential shape changes
along one side. (e) a piecewise-linear waveguide. (f) a waveguide with
a sinusoidal side variation. 
All waveguides except (b) are of length $1 \um$.
Their thicknesses are fixed at $0.22 \um$. The speedup of our method over the conventional RCWA for 
the same accuracy level (i.e., $10^{-4}$) is marked in yellow color. See \textbf{Data File 3-8} for underlying values.}
    \label{fig:simple_results}
\end{figure}
\subsection{Performance for Various Waveguides}
Next, we test the performance of our method for analyzing various photonic structures.
The shapes of these structures are described in \figref{simple_results}.
For each structure, the ground-truth scattering matrix is computed using the conventional RCWA with 
a high resolution in $z$-direction ($N=1024$).
Provided the ground-truth scattering matrix, we measure the performance-accuracy
curve for our method and the conventional RCWA: for our method, we progressively reduce the error 
threshold (i.e., $\alpha$ in \algref{adaptive_variant_rcwa}); and for each error threshold, we measure the
resulting scattering matrix error and the computational cost. This allows us to plot a curve showing how the 
accuracy changes over the computation time. Similary, for conventional RCWA, we progressively increase the 
$z$-direction resolution, and measure the scattering matrix error and computation time. 
All the timings are measured on a workstation with 8 Intel Xeon(R) E5-1620 CPUs running at 3.60GHz and an NVIDIA GeForce GTX 1080 GPU.

The resulting performance-accuracy curves are shown in \figref{simple_results}.
Given a fixed accuracy level (indicated by the green horizontal line in
\figref{simple_results}), our method is significantly faster. For example, to
compute the scattering matrix with an error around
$10^{-4}$, out method is at least $3.5\times$ faster than the conventional RCWA
for all test cases. In certain cases such as \figref{simple_results}-e, our
method is even an order of magnitude faster.  In all cases, our method
converges faster than the conventional RCWA, confirming that ours is a
higher-order method. This means that as the speedup of our method will become
more pronounced as the desired accuracy level increases.



We also note that when the conventional RCWA is used for analyzing a particular 
photonic structure, there exists no guideline to determine the $z$-direction resolution for a
certain accuracy level. As a result, in practice, one has to rely on multiple trials to choose
a proper resolution, and thus spend more time than what is reported
in \figref{simple_results}.


\begin{figure}[!htp]
    \vspace{-9mm}
    \centering
    \includegraphics[width=0.95\textwidth]{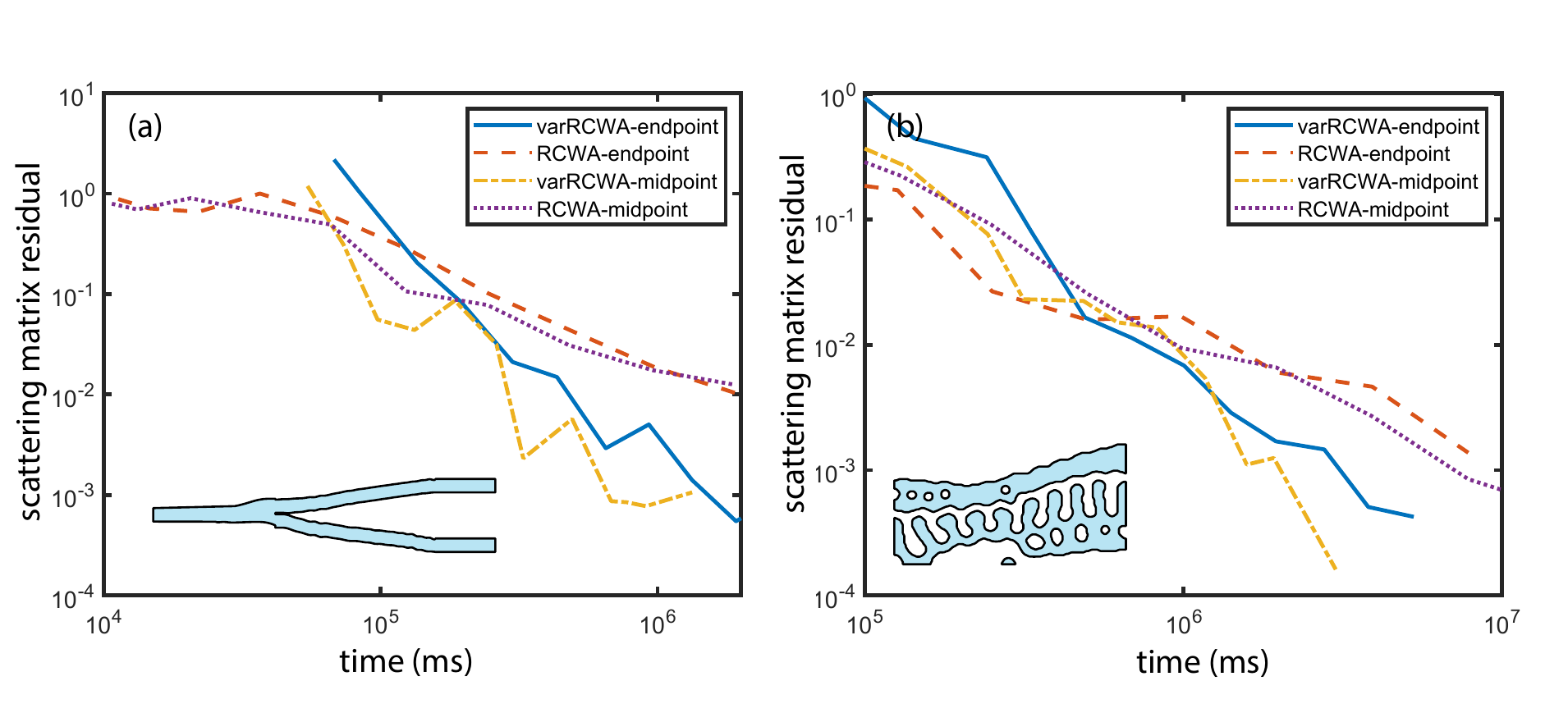}
    \vspace{-4mm}
    \caption{\textbf{Tests on complex photonic structures.} 
    (a) a Y-shape waveguide. (b) a frequency splitter reported in \cite{lumericalcontour}. See \textbf{Data File 9-12} for underlying values.}
    \label{fig:real_results}
\end{figure}

Finally, we evaluate our method on two real-world examples shown
in \figref{real_results}, namely a Y-shaped waveguide
and a frequency splitter. The former is smoothed out from our design of 90/10 splitter~\cite{zhu2021ultra} using B\'{e}zier curve and the latter is reported in ~\cite{lumericalcontour}.
Again, the ground-truth scattering matrices for both cases are obtained by
the conventional RCWA, with $8192$ discrete sections for the Y-shaped waveguide and
$32768$ sections for the frequency splitter, since the latter has more complex geometry.
As shown in \figref{real_results}, our method converges faster, 
and for sufficiently high accuracy, our method has much lower computational cost.


\section{Conclusion}
In summary, we present a high-order semi-analytical method.  Our method is best
suited for simulating photonic structures whose cross-sectional shapes vary
along the propagation direction. In comparison to conventional semi-analytical methods, 
our method is faster and more accurate.
In addition, our method allows the user to
specify an accuracy level, and adaptively discretizes the structure along $z$-direction
to achieve the desired accuracy.


In future, this work can be extended in various directions. Currently, when
subdivision of a section occurs, the section is subdivided into three subsections uniformly.
It is also possible to subdivide non-uniformly, according to the local cross-sectional shape 
of the structure.
Moreover, our implementation uses the first-order expansion that we derived,
although we can easily extend it to use even higher order expansion by recursively substituting 
a lower-order expansion into \eqref{eq:u_integral} and \eqref{eq:d_integral}.
In practice, our method based on the first-order expansion already outperforms the conventional 
RCWA.

\vspace{2mm}
\paragraph{Funding.}
National Science Foundation (1910839). 
\paragraph{Acknowledgments.}
We thank Utsav Dave and Michal Lipson for valuable suggestions.
\paragraph{Disclosures.}
The authors declare no conflicts of interest.

\vspace{2mm}
\noindent
See Supplement 1 for supporting content.
\bibliography{main}

\begin{thebibliography}{10}

\bibitem{blanes1998magnus}
{\sc Blanes, S., Casas, F., Oteo, J., and Ros, J.}
\newblock Magnus and fer expansions for matrix differential equations: the
  convergence problem.
\newblock {\em Journal of Physics A: Mathematical and General 31}, 1 (1998),
  259.

\bibitem{blanes2009magnus}
{\sc Blanes, S., Casas, F., Oteo, J.-A., and Ros, J.}
\newblock The magnus expansion and some of its applications.
\newblock {\em Physics reports 470}, 5-6 (2009), 151--238.

\bibitem{Chen14:ANM}
{\sc Chen, X., Zheng, C., Xu, W., and Zhou, K.}
\newblock An asymptotic numerical method for inverse elastic shape design.
\newblock {\em ACM Transactions on Graphics (Proceedings of SIGGRAPH 2014) 33},
  4 (Aug. 2014).

\bibitem{chu2003finite}
{\sc Chu, H.}
\newblock Finite difference approach to optical scattering of gratings.
\newblock In {\em Advanced Characterization Techniques for Optics,
  Semiconductors, and Nanotechnologies\/} (2003), vol.~5188, International
  Society for Optics and Photonics, pp.~358--370.

\bibitem{cochelin1994path}
{\sc Cochelin, B.}
\newblock A path-following technique via an asymptotic-numerical method.
\newblock {\em Computers \& structures 53}, 5 (1994), 1181--1192.

\bibitem{divitt2019ultrafast}
{\sc Divitt, S., Zhu, W., Zhang, C., Lezec, H.~J., and Agrawal, A.}
\newblock Ultrafast optical pulse shaping using dielectric metasurfaces.
\newblock {\em Science 364}, 6443 (2019), 890--894.

\bibitem{helton1976numerical}
{\sc Helton, J., and Stuckwisch, S.}
\newblock Numerical approximation of product integrals.
\newblock {\em Journal of Mathematical Analysis and Applications 56}, 2 (1976),
  410--437.

\bibitem{iserles1999solution}
{\sc Iserles, A., and N{\o}rsett, S.~P.}
\newblock On the solution of linear differential equations in lie groups.
\newblock {\em Philosophical Transactions of the Royal Society of London.
  Series A: Mathematical, Physical and Engineering Sciences 357}, 1754 (1999),
  983--1019.

\bibitem{jing2013analysis}
{\sc Jing, X., Jin, S., Tian, Y., Liang, P., Dong, Q., and Wang, L.}
\newblock Analysis of the sinusoidal nanopatterning grating structure.
\newblock {\em Optics \& Laser Technology 48\/} (2013), 160--166.

\bibitem{li2020efficient}
{\sc Li, J., Shi, L.~H., Ma, Y., Ran, Y., Liu, Y., and Wang, J.}
\newblock Efficient and stable implementation of rcwa for ultrathin multilayer
  gratings: T-matrix approach without solving eigenvalues.
\newblock {\em IEEE Antennas and Wireless Propagation Letters\/} (2020).

\bibitem{liu2012s4}
{\sc Liu, V., and Fan, S.}
\newblock S4: A free electromagnetic solver for layered periodic structures.
\newblock {\em Computer Physics Communications 183}, 10 (2012), 2233--2244.

\bibitem{lumericalcontour}
{\sc Lumerical}.
\newblock {GDS} pattern extraction for inverse designed devices using contours
  method, 2021.

\bibitem{miller2020large}
{\sc Miller, S.~A., Chang, Y.-C., Phare, C.~T., Shin, M.~C., Zadka, M.,
  Roberts, S.~P., Stern, B., Ji, X., Mohanty, A., Gordillo, O. A.~J., et~al.}
\newblock Large-scale optical phased array using a low-power multi-pass silicon
  photonic platform.
\newblock {\em Optica 7}, 1 (2020), 3--6.

\bibitem{mohamad2020fast}
{\sc Mohamad, H., Essaidi, S., Blaize, S., Macias, D., Benech, P., and Morand,
  A.}
\newblock Fast fourier factorization for differential method and rcwa: a
  powerful tool for the modeling of non-lamellar metallic diffraction gratings.
\newblock {\em Optical and Quantum Electronics 52}, 2 (2020), 1--13.

\bibitem{moharam1981rigorous}
{\sc Moharam, M., and Gaylord, T.}
\newblock Rigorous coupled-wave analysis of planar-grating diffraction.
\newblock {\em JOSA 71}, 7 (1981), 811--818.

\bibitem{cuda}
{\sc NVIDIA, Vingelmann, P., and Fitzek, F.~H.}
\newblock Cuda, release: 10.2.89, 2020.

\bibitem{piggott2015inverse}
{\sc Piggott, A.~Y., Lu, J., Lagoudakis, K.~G., Petykiewicz, J., Babinec,
  T.~M., and Vu{\v{c}}kovi{\'c}, J.}
\newblock Inverse design and demonstration of a compact and broadband on-chip
  wavelength demultiplexer.
\newblock {\em Nature Photonics 9}, 6 (2015), 374--377.

\bibitem{pregla1989method}
{\sc Pregla, R., and Pascher, W.}
\newblock The method of lines.
\newblock {\em Numerical techniques for microwave and millimeter wave passive
  structures 1\/} (1989), 381--446.

\bibitem{redheffer1959inequalities}
{\sc Redheffer, R.}
\newblock Inequalities for a matrix riccati equation.
\newblock {\em Journal of Mathematics and Mechanics\/} (1959), 349--367.

\bibitem{roberts2017rigorous}
{\sc Roberts, C.~M., and Podolskiy, V.~A.}
\newblock Rigorous diffraction interface theory.
\newblock {\em Applied Physics Letters 110}, 17 (2017), 171108.

\bibitem{slavik2007product}
{\sc Slav{\'\i}k, A.}
\newblock {\em Product integration, its history and applications}.
\newblock Matfyzpress Prague, 2007.

\bibitem{takegoshi2015comparison}
{\sc Takegoshi, K., Miyazawa, N., Sharma, K., and Madhu, P.}
\newblock Comparison among magnus/floquet/fer expansion schemes in solid-state
  nmr.
\newblock {\em The Journal of chemical physics 142}, 13 (2015), 134201.

\bibitem{yee1966numerical}
{\sc Yee, K.}
\newblock Numerical solution of initial boundary value problems involving
  maxwell's equations in isotropic media.
\newblock {\em IEEE Transactions on antennas and propagation 14}, 3 (1966),
  302--307.

\bibitem{zhu2020inverse}
{\sc Zhu, Z., Dave, U.~D., Lipson, M., and Zheng, C.}
\newblock Inverse geometric design of fabrication-robust nanophotonic
  waveguides.
\newblock In {\em 2020 Conference on Lasers and Electro-Optics (CLEO)\/}
  (2020), IEEE, pp.~1--2.

\bibitem{zhu2021ultra}
{\sc Zhu, Z., Dave, U.~D., Lipson, M., and Zheng, C.}
\newblock Ultra-broadband nanophotonics via adaptive inverse design.
\newblock In {\em 2021 Conference on Lasers and Electro-Optics (CLEO)\/}
  (2021), IEEE, pp.~1--2.

\end{thebibliography}


\begin{thebibliography}{1}

\bibitem{iserles2004method}
{\sc Iserles, A.}
\newblock On the method of neumann series for highly oscillatory equations.
\newblock {\em Bit Numerical Mathematics 44}, 3 (2004), 473--488.

\end{thebibliography}
\bibliographystyle{acm}
\end{document}


\maketitle
\section{Explanation of the requirement of the integrand in Eqn.~(6)}
To use Eqn.~(6) for a long section, the coefficients
$\mat{W}^{-1}\mat{P}\mat{V}$ and $\mat{V}^{-1}\mat{Q}\mat{W}$ should be small
enough. We explain the intuition behind this requirement using a simple
example, a scalar differential equation in a similar form of Eqn.~(5)
\begin{equation}\label{eqn:simple_diff}
    \frac{\text{d}}{\text{d}z}y(z) = a y(z),
\end{equation}
where $a$ is a constant, similar to 
$\mat{W}^{-1}\mat{P}\mat{V}$ and $\mat{V}^{-1}\mat{Q}\mat{W}$ in
 Eqn.~(5),
and $y(z)$ is
similar to the vectors $\bvec{a}(z')$ and $\bvec{b}(z')$ therein. 
This equation can be written as an integral form similar to Eqn.~(6):
 \begin{equation}\label{eqn:simple_integral}
     y(z) = y(z_L) + \int^z_{z_L}ay(z')\diff z'.
 \end{equation}
By recursively substituting \eqref{eqn:simple_integral} into $y(z')$, we get
 \begin{equation}
    y(z) =\left[1 +  \int a\diff z'  + \frac{1}{2}\left(\int a\diff z'\right)^2 + \cdots\right] y(z_L),
 \end{equation}
where the right-hand side is a Taylor expansion of $\exp{\left(\int a\diff z'\right)}$:
 \begin{equation}
     \exp{\left(\int a\diff z'\right)} = 1 +  \int a\diff z'  + \frac{1}{2}\left(\int a\diff z'\right)^2 + \cdots
 \end{equation}
Although Taylor expansion of an exponential always converges for all $\int
a\diff z'$, if we truncate this series to the first order, the residual, which is
\begin{equation}
    \text{Re} = \left|\exp{\left(\int a\diff z'\right)} - 1 -  \int a\diff z'\right|
    = \left|\left(\int a\diff z'\right)^2\right|\left|\frac{1}{2} + o\left(\int a\diff z'\right)\right|,
\end{equation}
scales with the magnitude of $\int a\diff z' = a(z - z_L) = aL$. In order to
increase the length $L$ of this integral while maintaining the same integral
value $aL$ (hence the accuracy), $a$ should be small enough. Similarly
$\mat{W}^{-1}\mat{P}\mat{V}$ and $\mat{V}^{-1}\mat{Q}\mat{W}$ should be small
enough to ensure plausible accuracy of the first order truncation. $\mat{W}$ and
$\mat{V}$ do not scale with the section length as they are eigenvectors.
Therefore, $\mat{P}$ and $\mat{Q}$ must have small enough norms, 
which unfortunately can not be guaranteed.

\section{Derivation of Eqn.~(8)}
Note that the coefficients $\mat{P}$ and $\mat{Q}$, which depend on the cross-sectional material distribution,
can be arbitrarily large, even if
the section has a fixed cross section. 
We therefore wish to replace them by $\delta
\mat{P}$ and $\delta \mat{Q}$, which depend only on the cross-sectional variation of the section.
To this end, instead of differentiating $\bvec{a}$ and $\bvec{b}$, we scale
$\bvec{a}$ and $\bvec{b}$ first in order to cancel out some terms and construct $\delta
\mat{P}$ and $\delta \mat{Q}$. In particular, let
\begin{equation}
    \tilde{\bvec{a}} = \gamma_1(z)\bvec{a}(z)
\end{equation}
\begin{equation}
    \tilde{\bvec{b}} = \gamma_2(z)\bvec{b}(z)
\end{equation}
for some specific choices of $\gamma_1(z)$ and $\gamma_2(z)$.
Inspired by \cite{iserles2004method} and the derivation of RCWA, after several
trials, we find $\gamma_1(z)$ and $\gamma_2(z)$ in Eqs.~(7).

Firstly, we differentiate $\tilde{\bvec{a}}$ in the first part of Eqn.~(7) as
\begin{equation}\label{eq:differentiation_u}
    \frac{\partial \tilde{\bvec{a}}}{\partial z} = -\frac{j}{k_0}\mat{\Lambda}\tilde{\bvec{a}} + \exp{\left(-\frac{j}{k_0}\mat{\Lambda}z\right)}\frac{\partial \bvec{a}}{\partial z}.
\end{equation}
The first term $-\frac{j}{k_0}\mat{\Lambda}\tilde{\bvec{a}}$ is emerged to 
construct $\delta \mat{P}$ and $\delta \mat{Q}$.
Now substitude Eqs.~(5) in \eqref{eq:differentiation_u}: 
\begin{equation}\label{eq:differentiation_u_2}
    \frac{\partial \tilde{\bvec{a}}}{\partial z} = \frac{j}{2k_0}\exp{\left(-\frac{j}{k_0}\mat{\Lambda}z\right)}\left(-2\mat{\Lambda}\bvec{a} + \mat{W}^{-1}\mat{P}\mat{V}\left(\bvec{a} - \bvec{b}\right) + \mat{V}^{-1}\mat{Q}\mat{W}\left(\bvec{a} + \bvec{b}\right)\right).
\end{equation}
From $\mat{P}_r\mat{Q}_r = \mat{W}\mat{\Lambda}^2\mat{W}^{-1}$ and $\mat{V} = \mat{Q}_r\mat{W}\mat{\Lambda}^{-1}$, we get the following relations:
\begin{equation}\label{eq:lambda_in_P_Q}
    \mat{\Lambda} = \mat{V}^{-1}\mat{Q}_r\mat{W} = \mat{W}^{-1}\mat{P}_r\mat{V}.
\end{equation}

Next, by substituting Eqs.~(\ref{eq:lambda_in_P_Q}) into \eqref{eq:differentiation_u_2}, we obtain
\begin{equation*}\label{eq:differentiation_u_3}
    \begin{split}
    \frac{\partial \tilde{\bvec{a}}}{\partial z} 
    &= \frac{j}{2k_0}\exp{\left(-\frac{j}{k_0}\mat{\Lambda}z\right)}\left(\mat{W}^{-1}\left(\mat{P}-\mat{P}_r\right)\mat{V}\bvec{a} + \mat{V}^{-1}\left(\mat{Q}-\mat{Q}_r\right)\mat{W}\bvec{a}-\mat{W}^{-1}\mat{P}\mat{V}\bvec{b} + \mat{V}^{-1}\mat{Q}\mat{W}\bvec{b}\right)\\
    &= \frac{j}{2k_0}\exp{\left(-\frac{j}{k_0}\mat{\Lambda}z\right)}\left(\mat{W}^{-1}\delta\mat{P}\mat{V}\bvec{a} + \mat{V}^{-1}\delta\mat{Q}\mat{W}\bvec{a}-\mat{W}^{-1}\mat{P}\mat{V}\bvec{b} + \mat{V}^{-1}\mat{Q}\mat{W}\bvec{b}\right)
    \end{split}
\end{equation*}
There are still $\mat{P}$ and $\mat{Q}$ in the coefficients of $\bvec{b}$.
But the signs of the coefficients are different, one is positive and the other
one is negative. Therefore, we can add an extra zero term $(\mat{\Lambda} -
\mat{\Lambda})\bvec{b}$ to help with the cancellation:
\begin{equation}\label{eq:differentiation_u_4}
  \begin{split}
    \frac{\partial \tilde{\bvec{a}}}{\partial z} 
    &= \frac{j}{2k_0}\exp{\left(-\frac{j}{k_0}\mat{\Lambda}z\right)}\left(\delta\mat{A}\bvec{a}-\mat{W}^{-1}\mat{P}\mat{V}\bvec{b} + \mat{V}^{-1}\mat{Q}\mat{W}\bvec{b} + \left(\mat{\Lambda} - \mat{\Lambda}\right)\bvec{b}\right)\\
    &= \frac{j}{2k_0}\exp{\left(-\frac{j}{k_0}\mat{\Lambda}z\right)}\left(\delta\mat{A}\bvec{a}-\mat{W}^{-1}\mat{P}\mat{V}\bvec{b} + \mat{V}^{-1}\mat{Q}\mat{W}\bvec{b} + \left(\mat{W}^{-1}\mat{P}_r\mat{V}-  \mat{V}^{-1}\mat{Q}_r\mat{W}\right)\bvec{b}\right)\\
    &= \frac{j}{2k_0}\exp{\left(-\frac{j}{k_0}\mat{\Lambda}z\right)}\left(\delta\mat{A}\bvec{a}-\delta\mat{B}\bvec{b}\right)\\
    &= \frac{j}{2k_0}\exp{\left(-\frac{j}{k_0}\mat{\Lambda}z\right)}\delta\mat{A}\exp{\left(\frac{j}{k_0}\mat{\Lambda}z\right)}\tilde{\bvec{a}}-\frac{j}{2k_0}\exp{\left(-\frac{j}{k_0}\mat{\Lambda}z\right)}\delta\mat{B}\exp{\left(-\frac{j}{k_0}\mat{\Lambda}z\right)}\tilde{\bvec{b}}.
  \end{split}
\end{equation}
Thanks to the symmetry of $\tilde{\bvec{a}}$ and $\tilde{\bvec{b}}$, we can get the differentiation of $\tilde{\bvec{b}}$ the same way.

\section{Construction of the high order scattering matrix}
The matrix elements in Eqn.(3) can be extracted from Eqs.~(12) and (13) as follows:
\begin{align}
    \mat{T}_{LR} &= e^{\frac{j}{k_0}\mat{\Lambda}(z_R-z_L)} + \frac{j}{2k_0}\int^{z_R}_{z_L}e^{\frac{j}{k_0}\mat{\Lambda}(z_R-z')}\delta\mat{A}(z')e^{\frac{j}{k_0}\mat{\Lambda}(z'-z_L)}\diff z' \label{eq:smat_tlr} \\
    \mat{R}_R &= - \frac{j}{2k_0}\int^{z_R}_{z_{L}}e^{\frac{j}{k_0}\mat{\Lambda}(z_R-z')}\delta \mat{B}(z')e^{\frac{j}{k_0}\mat{\Lambda}(z_R-z')}\diff z' \label{eq:smat_rr} \\
    \mat{R}_L &= -\frac{j}{2k_0}\int^{z_R}_{z_L}e^{\frac{j}{k_0}\mat{\Lambda}(z'-z_L)}\delta\mat{B}(z')e^{\frac{j}{k_0}\mat{\Lambda}(z'-z_L)}\diff z' \label{eq:smat_rl} \\
    \mat{T}_{RL} &=  e^{\frac{j}{k_0}\mat{\Lambda}(z_R-z_L)} + \frac{j}{2k_0}\int^{z_R}_{z_L}e^{\frac{j}{k_0}\mat{\Lambda}(z'-z_L)}\delta \mat{A}(z')e^{\frac{j}{k_0}\mat{\Lambda}(z_R-z')}\diff z'. \label{eq:smat_trl}
\end{align}

\section{Projection of scattering matrix}
Because the scattering matrix at a section $i$ in (given in Eqs.~(\ref{eq:smat_tlr}-\ref{eq:smat_trl}))
is expressed in the basis of eigenmodes at the reference position of section $i$.
If the eigenmodes are different from those of section $(i-1)$, we need to
reproject the scattering matrix at section $i$ to match the eigenmodes of
section $(i-1)$ before using Redheffer star product to multiply them
together.

According to the continuity of the tangential fields at the interface between the two sections, if we denote the eigenmodes at section $i$ as $\mat{W}_i$ and $\mat{V}_i$, and we use the superscript $(i-1)$ to denote the vector is represented in the basis of section $(i-1)$, while the vectors without superscripts are in the basis of the section $i$. From Eq.~(2), we have
\begin{equation}\label{eq:WV_eq}
    \begin{split}
     \mat{W}_{i-1}\left(\bvec{a}^{(i-1)}_L + \bvec{b}^{(i-1)}_L\right) &= \mat{W}_i\left(\bvec{a}_L + \bvec{b}_L\right) \\
     \mat{V}_{i-1}\left(\bvec{a}^{(i-1)}_L - \bvec{b}^{(i-1)}_L\right) &= \mat{V}_i\left(\bvec{a}_L - \bvec{b}_L\right).
    \end{split}
\end{equation}
Now we denote 
\begin{equation}
    \begin{split}
     \mat{X} &= \frac{1}{2}\left(\mat{W}^{-1}_i\mat{W}_{i-1}+\mat{V}^{-1}_i\mat{V}_{i-1}\right)\\
     \mat{Y} &= \frac{1}{2}\left(\mat{W}^{-1}_i\mat{W}_{i-1}-\mat{V}^{-1}_i\mat{V}_{i-1}\right),
\end{split}
\end{equation}
which allow us to express $\bvec{a}_L$ and $\bvec{b}_L$ using \eqref{eq:WV_eq}, namely,
\begin{equation}\label{eq:s18}
    \begin{split}
     \bvec{a}_L &= \mat{X}\bvec{a}^{(i-1)}_L+ \mat{Y}\bvec{b}^{(i-1)}_L \\
    \bvec{b}_L  &= \mat{Y}\bvec{a}^{(i-1)}_L+ \mat{X}\bvec{b}^{(i-1)}_L.
    \end{split}
\end{equation}
Meanwhile, according to the scattering matrix definition in Eq.~(3), we have
\begin{equation}\label{eq:s19}
    \begin{split}
    \bvec{b}_L &= \mat{R}_L\bvec{a}_L + \mat{T}_{RL}\bvec{b}_R \\
 \left(\mat{X}- \mat{R}_L\mat{Y}\right)\bvec{b}^{(i-1)}_L &=-\left(\mat{Y} - \mat{R}_L\mat{X}\right)\bvec{a}^{(i-1)}_L+ \mat{T}_{RL}\bvec{b}_R.
    \end{split}
\end{equation}
By comparing \eqref{eq:s18} to \eqref{eq:s19}, we obtain the recipe for reprojection of the scattering matrix:
\begin{equation}
\begin{split}
 \tilde{\mat{R}}_{L} &=-\left(\mat{X}- \mat{R}_L\mat{Y}\right)^{-1}\left(\mat{Y} - \mat{R}_L\mat{X}\right)\\
 \tilde{\mat{T}}_{RL}&= \left(\mat{X}- \mat{R}_{L}\mat{Y}\right)^{-1}\mat{T}_{RL}.
\end{split}
\end{equation}
Similarly, from
\begin{equation}
\bvec{a}_R=\mat{T}_{LR}\mat{X}\bvec{a}^{(i-1)}_L+  \mat{T}_{LR}\mat{Y}\left(\tilde{\mat{R}}_{L}\bvec{a}^{(i-1)}_L +\tilde{\mat{T}}_{RL}\bvec{b}_R \right)+\mat{R}_{R}\bvec{b}_R,
\end{equation}
We have
\begin{equation}
\begin{split}
 \tilde{\mat{T}}_{LR} &=\mat{T}_{LR}\mat{X}+\mat{T}_{LR}\mat{Y}\tilde{\mat{R}}_{L} \\
 \tilde{\mat{R}}_R    &= \mat{T}_{LR}\mat{Y}\tilde{\mat{T}}_{RL} + \mat{R}_{R}.
\end{split}
\end{equation}
The projected scattering matrix elements $\tilde{\mat{R}}_{L}$, $\tilde{\mat{T}}_{RL}$, $\tilde{\mat{T}}_{LR}$, $\tilde{\mat{R}}_{R}$ can be left multiplied by the scattering matrix of section $(i-1)$.
\section{CPU Implementation}
\begin{figure}[!htp]
    \centering
    \includegraphics[width=0.7\textwidth]{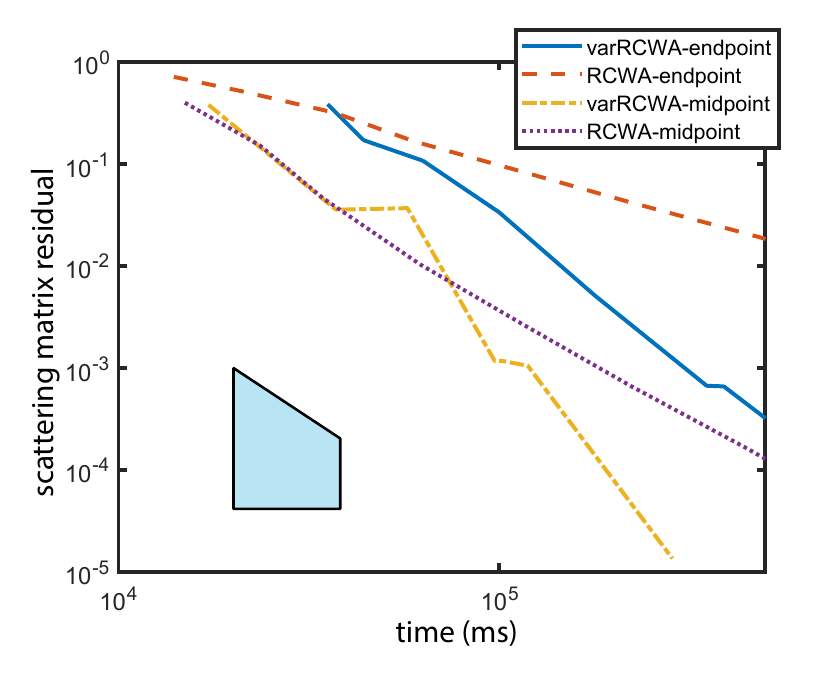}
    \vspace{-5mm}
    \caption{CPU implementation of the experiments in Fig.~5(a). See \textbf{Data File 13} for underlying values.}
    \label{fig:cpu_implementation}
\end{figure}
We have also implemented our method on CPU. We choose the same structure as
Fig.~5-a to run the experiments. As is shown in
Fig~\ref{fig:cpu_implementation}, our method performs better no matter we
choose the midpoint or endpoint as the reference point. However, because the
program runs sequentially, the CPU implementation cannot parallelize the
computation of matrix products, which is not so efficient. Therefore, matrix
product still takes a noticeable amount of time ($1/20$ per operation) compared
to eigenvalue decomposition. Although both RCWA and our method running on CPU
are slower, our method is still faster compared to conventional RCWA, although
the speedup of our method is not as significant as on GPU. 

\section{Error Estimation}
When subdividing the waveguide, our algorithm checks if an estimated error meets
the user-provided accuracy level $\alpha$.
To estimate the error, we compare
the scattering matrix of the $p$th order and the $(p-1)$th order. In our case, we subtract the results standard RCWA ($p=0$) from the results in Eqs.~(\ref{eq:smat_tlr}-\ref{eq:smat_trl}) ($p=1$), and get
    the following four difference matrices
\begin{align}
    \Delta\mat{T}_{LR} &=  \frac{j}{2k_0}\int^{z_R}_{z_L}e^{\frac{j}{k_0}\mat{\Lambda}(z_R-z')}\delta\mat{A}(z')e^{\frac{j}{k_0}\mat{\Lambda}(z'-z_L)}\diff z' \label{eq:smat_tlr} \\
    \Delta\mat{R}_R &= - \frac{j}{2k_0}\int^{z_R}_{z_{L}}e^{\frac{j}{k_0}\mat{\Lambda}(z_R-z')}\delta \mat{B}(z')e^{\frac{j}{k_0}\mat{\Lambda}(z_R-z')}\diff z' \label{eq:smat_rr} \\
    \Delta\mat{R}_L &= -\frac{j}{2k_0}\int^{z_R}_{z_L}e^{\frac{j}{k_0}\mat{\Lambda}(z'-z_L)}\delta\mat{B}(z')e^{\frac{j}{k_0}\mat{\Lambda}(z'-z_L)}\diff z' \label{eq:smat_rl} \\
    \Delta\mat{T}_{RL} &= \frac{j}{2k_0}\int^{z_R}_{z_L}e^{\frac{j}{k_0}\mat{\Lambda}(z'-z_L)}\delta \mat{A}(z')e^{\frac{j}{k_0}\mat{\Lambda}(z_R-z')}\diff z'. \label{eq:smat_trl}
\end{align}
We use the maximum value $\varepsilon$ of the max norms of the above four matrices as the estimated error, and compare $\varepsilon$ with user-provided $\alpha$.

\begin{figure}[!htp]
    \centering
    \vspace{-50mm}
    \includegraphics[width=\textwidth]{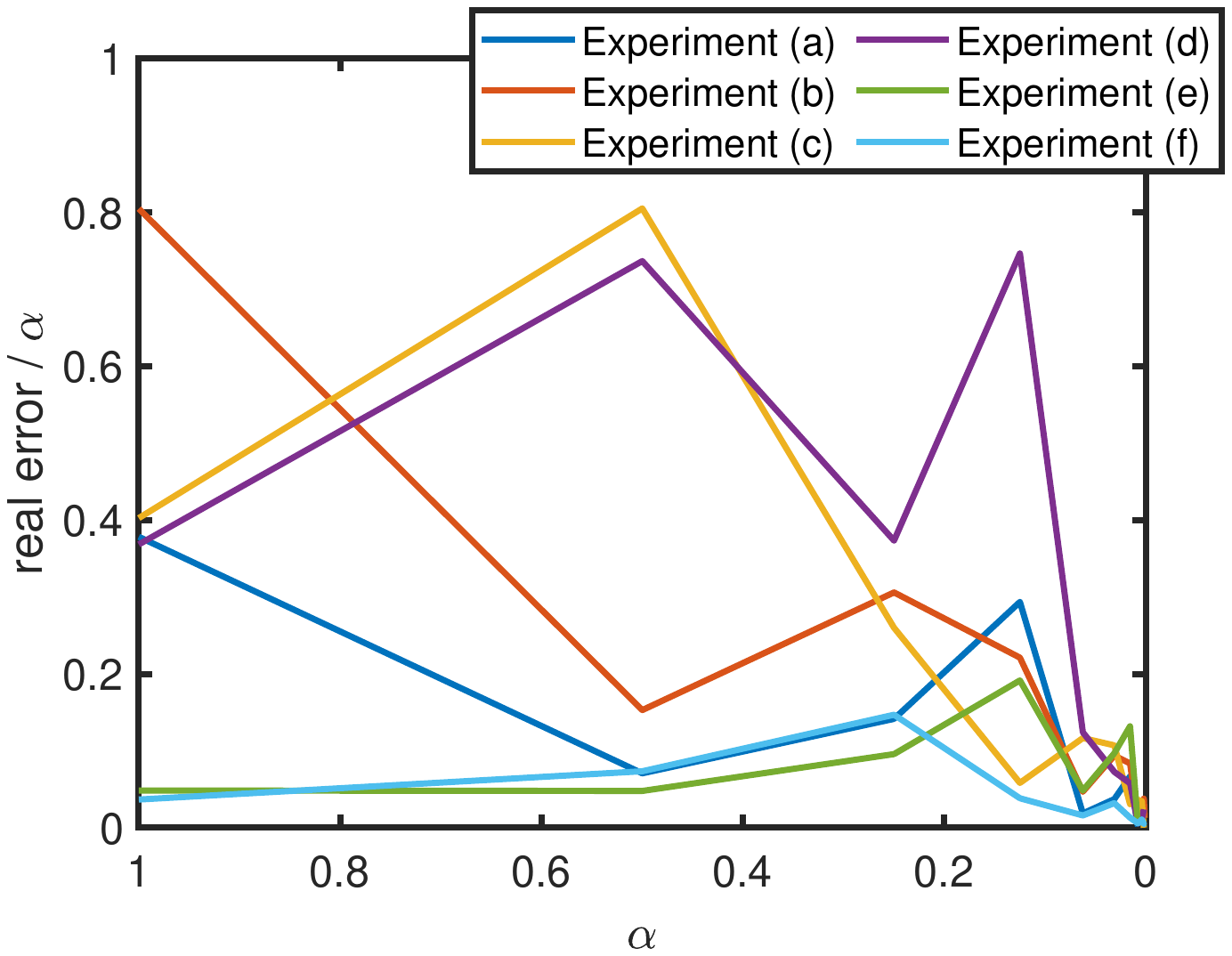}
    \vspace{-55mm}
    \caption{\textbf{Real errors versus different $\alpha$s for experiments in Fig. 5.} We rerun each experiment with different user-provided $\alpha$s ranging from $0$ to $1$. When the algorithm terminates, we compare each scattering matrix with the corresponding ground truth scattering matrix (running at a large number of sections $N=1024$) to get the \textit{real error} (max norm of the difference between two matrices). It shows the real error of the scattering matrix is always bounded by $\alpha$. See \textbf{Data File 14} for underlying values.}
    \label{fig:alpha_vs_real_erro}
\end{figure}

Although $\varepsilon$ is estimated, any error caused by higher order perturbations should be smaller than  $\varepsilon$. 
We have also compared this estimated error with the final
error in the scattering matrix for each experiments in Fig.~5. As is shown in
Fig.~\ref{fig:alpha_vs_real_erro}, for each user-specified $\alpha$, the real error is always smaller
than $\alpha$, which indicates it is a good overestimate.

\bibliography{main}
\bibliographystyle{acm}